\documentclass[11pt]{article}

\def\init{\setcounter{equation}{0}}

\newtheorem{theoreme}{Theorem}[section]
\newtheorem{proposition}[theoreme]{Proposition}
\newtheorem{lemme}[theoreme]{Lemma}
\newtheorem{definition}[theoreme]{Definition}
\newtheorem{corollaire}[theoreme]{Corollary}
\newtheorem{conjecture}[theoreme]{Conjecture}

\newcommand{\bsq}{{\vrule height .9ex width .8ex depth -.1ex }}
\def\rr{{\bf R}}

\def\nn{{\bf N}}
\def\zz{{\bf Z}}
\def\qq{{\bf Q}}

\def\ie{{\it i.e., }}
\def\eg{{\it e.g., }}

\def\proof{{\it  Proof. }}
\def\qed{$\bsq$}

\def\half{\frac{1}{2}}
\def\eps{\epsilon}

\def\coinf{C_{0}^{\infty}}

\def\mod{\hbox{mod}\,}

\def\deg{\hbox{deg}\,}

\newcommand{\beq}{\begin{equation}}
\newcommand{\eeq}{\end{equation}}
\newcommand{\bet}{\begin{theoreme}}
\newcommand{\eet}{\end{theoreme}}
\newcommand{\bear}[1]{\begin{array}{#1}}
\newcommand{\ear}{\end{array}}

\def\aa{\alpha}
\def\l{\lambda}
\def\eps{\epsilon}

\begin{document}

\title{The spectral set conjecture and multiplicative properties
of roots of polynomials}
\author{I. \L aba
\thanks{Supported in part by NSERC}
\\ Department of Mathematics\\University of British Columbia\\
Vancouver, Canada V6T 1Z2\\ {\it ilaba@math.ubc.ca}}
\maketitle

\begin{abstract}
Fuglede's conjecture \cite{Fug} states that a set $\Omega\subset
\rr^n$ tiles $\rr$ by translations if and only if $L^2(\Omega)$
has an orthogonal basis of exponentials.  We obtain new partial
results supporting the conjecture in dimension 1.
\end{abstract}

\section{The results}
\label{intro}
\init

A Borel set $\Omega\subset\rr^n$ of positive measure is said to 
{\it tile $\rr^n$ by translations} if there is a discrete set 
$T\subset\rr^n$ such that, up to sets of measure 0, the sets $\Omega
+t,\ t\in T,$ are disjoint and $\bigcup_{t\in T}(\Omega+t)=\rr^n$.
We will refer to $T$ as the {\em translation set}, and write $\rr=\Omega
\oplus T$.
We also say that $\Lambda=\{\lambda_k\}_{k\in\zz}\subset\rr^n$ is
a {\it spectrum} for $\Omega$ if:
\beq
\{e^{2\pi i\l_k\cdot x}\}_{k\in\zz}\hbox{ is an orthogonal basis for }
L^2(\Omega).
\label{a.e00}
\eeq
A {\it spectral set} is a domain $\Omega\subset\rr^n$ such that (\ref{a.e00})
holds for some $\Lambda$.  
 
The following conjecture is due to Fuglede \cite{Fug}.

\begin{conjecture} {\em (The spectral set conjecture \cite{Fug}.)} A domain
$\Omega\subset\rr^n$ is a spectral set if and only if it tiles $\rr^n$
by translations.
\label{fug-conj}
\end{conjecture}

Conjecture \ref{fug-conj} has a functional-analytic origin: it
was motivated by a question of I. Segal concerning the
``extension property", \ie the existence of commuting self-adjoint extensions
of the operators $-i\frac{\partial}{\partial x_j}$, $j=1,\dots,n$,
defined on $\coinf(\Omega)$, to a dense subspace of $L^2(\Omega)$.
It turns out that any spectral set $\Omega$ has the extension
property; moreover, if $\Omega$ is assumed to be connected,
it must be a spectral set in order for the 
extension property to hold (\cite{Fug}, \cite{J}, \cite{P1}).  

Fuglede proved in \cite{Fug} that the conjecture is true if 
either $\Lambda$ or $T$ is a group. 
Subsequent attempts to resolve 
the general case have revealed connections to
functional analysis, representation theory, combinatorics, 
commutative algebra, and Fourier analysis, among others.
Fuglede's conjecture has also led to many other questions of
independent interest concerning the relations between the tiling and
spectral properties of sets, some of which have now been investigated
in detail.
(See \eg \cite{IKP}, \cite{IKT1}, \cite{IP}, \cite{JP1}, \cite{JP2},
\cite{JP4}, \cite{K1}, \cite{K2}, \cite{LRW}, \cite{LS}, \cite{LW1}, \cite{LW2}.)
 
For convex domains $\Omega\subset\rr^n$, the problem is now understood
to be related to the geometry of the set 
$\{\xi: \ \hat\chi_\Omega(\xi)=0\}$ (\cite{K1}, \cite{IKP}, \cite{IKT1}),
and the 2-dimensional convex case appears to be resolved \cite{IKT2}.
The general case is much more complicated, even in dimension 1, and
is still nowhere near resolution.  

The purpose of this paper is to relate the spectral set conjecture
for domains $\Omega=A+[0,1]\subset\rr$ to purely algebraic questions concerning
the multiplicative properties of roots of certain types of
polynomials, and to provide partial answers to these questions.

We briefly summarize some of the previous work relevant to the subject.
Lagarias and Wang \cite{LW1} proved that if $\Omega\subset\rr$ is a bounded
domain and $\Omega\oplus T=\rr$, then $T$ is rational and periodic.
The question of when a given periodic set $\Lambda$ is a spectrum
for a given tile set $T$ was discussed in \cite{LW2}, \cite{P2}.
In particular, it is proved in \cite{LW2} (Theorem 1.2) that if
$\Omega\oplus T=\rr$, $T=B+N\zz$, $B\subset\zz$, and the cyclic 
group $\zz_N$ has the {\em strong Tijdeman property}, then 
$\Omega$ is a spectral set; however, there are examples of $\zz_N$ 
which do not have the Tijdeman property, see \cite{Sz}, \cite{LS}.
Our Corollary \ref{cor1} uses the recent results of \cite{CM} to
improve on Theorem 1.3 of \cite{LW2}.  

It is not known whether a spectrum $\Lambda$ of a bounded domain
$\Omega\subset\rr$ must always be rational and periodic.
For domains of the form $\Omega=A+[0,1)$, the results of \cite{JP3}
(see Proposition \ref{JP-prop} below) imply that $\Lambda$ must indeed be periodic.
The rationality of $\Lambda$ seems to be a more difficult question:
Theorem \ref{thm-2N} below is a modest partial result in this direction.

Several special cases of the conjecture have
been investigated in more detail, see \eg \cite{L1}, \cite{PW}.
Note that the articles \cite{LW2}, \cite{PW} considered also
the related question of the existence of a {\em universal spectrum},
\ie a common spectrum for all domains $\Omega$ which tile $\rr^n$
by the same set of translations $T$. This question will not
be addressed in the present paper.

Let 
\beq
\Omega=\bigcup_{j=0}^{N-1}[a_j,a_j+1)\subset\rr,\ 
\Lambda=\{\l_k\}_{k\in\zz}\subset\rr.
\label{e.omega}
\eeq
We will say that $\Lambda$
is a spectrum for $\Omega$ if (\ref{a.e00}) holds.
Clearly, we may assume that $a_0=0$, $\l_0=0$. 

\begin{definition}
Let $A(x)$ be a polynomial of the form $A(x)=\sum_{k=0}^{N-1}
x^{a_k}$, where $a_k$ are distinct non-negative integers.
We will say that $\{\theta_1,\dots,\theta_{N-1}\}\subset(0,1)$
is a {\em spectrum } for $A(x)$ if the $\theta_j$ are all
distinct and:
\beq
A(\eps_{ij})=0\hbox{ for all }0\leq i,j\leq N-1,\ i\neq j,
\label{spec-pol}
\eeq
where 
\beq
\eps_{ij}=e^{2\pi i(\theta_i-\theta_j)},\ \theta_0=1.
\label{spec-pol-bis}
\eeq
\end{definition}

The reason for this terminology is the following result,
due to Jorgensen and Pedersen \cite{JP3} (see also \cite{P2}).

\begin{proposition}
(Jorgensen - Pedersen \cite{JP3})
Let $\Omega$ and $\Lambda$ be as in (\ref{e.omega}), $0\in\Lambda$.
Then $\Lambda$ is a spectrum 
for $\Omega$ if and only if all of the following are satisfied:

\smallskip
(i) $a_0,a_1,\dots,a_{N-1}\in\zz$;

\smallskip
(ii) the polynomial $A(x)=\sum_{k=0}^{N-1} x^{a_k}$
has a spectrum $\{\theta_j:\ j=1,\dots,N-1\}\subset(0,1)$;

\smallskip
(iii) $\Lambda=\bigcup_{j=0}^{N-1}(\theta_j+\zz)$, where
$\theta_0=0$ and $\theta_j$, $j=1,\dots,N-1$, are as in
{\it(ii)}.

\label{JP-prop}
\end{proposition}

This result was stated in \cite{JP3} in terms of Hadamard matrices:
it is easy to see that the condition (ii) above is
satisfied if and only if the columns of the matrix $(e^{2\pi ia_j\theta_k}
)_{j,k}$ are mutually orthogonal, which was the condition given in 
\cite{JP3}. Our reformulation of it in terms of
polynomials was motivated
by the recent work of Coven and Meyerowitz \cite{CM},
who related the tiling properties of a set $A\subset\zz$
to the algebraic properties of the corresponding 
polynomial $A(x)=\sum_{a\in A}x^a$.  (Note that 
such polynomials were also used in \cite{PW}.)
Clearly, $A\subset\zz$ tiles $\zz$ if and only if
$A+[0,1)$ tiles $\rr$, hence the relevance of this work
to the problem under consideration.  

The main result of \cite{CM} is as follows.
Recall that, for $s\in\zz$, the $s$-th cyclotomic polynomial
$\Phi_s(x)$ is defined inductively by
$x^s-1=\prod_{k|s}\Phi_s(x)$; 
equivalently, $\Phi_s(x)=\prod (x-\eps_i)$, where $\eps_i$
are the $s$-th primitive roots of 1.
Define $A(x)$ as above, and let $S_A$ be the set of prime
powers $s$ such that $\Phi_s(x)$ divides $A(x)$. 
Consider the following conditions on $A(x)$:

\medskip
(T1) $A(1)=\prod_{s\in S_A}\Phi_s(1)$;

\smallskip
(T2) If $s_1,\dots,s_k\in S_A$ are powers of different
primes, then $\Phi_{s_1\dots s_k}(x)$ divides $A(x)$.

\begin{theoreme}
(Coven - Meyerowitz \cite{CM})

\smallskip
(i) if (T1), (T2) hold, then $A$ tiles $\zz$ by translations;

\smallskip
(ii) if $A$ tiles $\zz$ by translations, then (T1) holds;

\smallskip
(iii) if $A$ tiles $\zz$ by translations and $N=\#A$ has at most two
prime factors, then (T2) holds.
\label{CM-thm}
\end{theoreme}

It is not known whether (T2) is satisfied for all sets $A$ 
which tile $\zz$ by translations.  It does hold in all cases
known to the author, and, in particular, it holds for the
sets constructed by Szab\'o \cite{Sz} which were used to
disprove Tijdeman's conjecture.

Our first result is that (T1), (T2) also guarantee the
existence of a spectrum, and that a weak partial converse holds.
Here and in the sequel, $A$ is a
subset of $\zz$ with $0\in A$, $N=\#A\geq 2$, $A(x)$,
$S_A$ are defined as above, and $\Omega=A+[0,1)$. 

\begin{theoreme}
(i) If $A(x)$ satisfies (T1), (T2), it has a spectrum.

\smallskip
(ii) If $A(x)$ has a spectrum $\{\theta_1,\dots,\theta_{N-1}\}
\subset p^{-\alpha}\zz$, where $p$ is a prime, $\alpha\in\nn$,
then $N$ is a power of $p$, and (T1) holds.
\label{tt-thm}
\end{theoreme}

Combining Theorem \ref{tt-thm} with Proposition \ref{JP-prop}
and Theorem \ref{CM-thm}, we obtain the
following immediate corollary, which improves on Theorem
1.3 of \cite{LW2}. (Note that if $N=p^\alpha$ is a prime power, then
(T2) automatically holds.)

\begin{corollaire}
(i) If $\Omega$ tiles $\zz$ and $N$ has at most 2 prime factors, then
$\Omega$ is a spectral set.

\smallskip
(ii) If $\Omega$ is a spectral set with a spectrum $\Lambda\subset
s^{-1}\zz$, where $s=p^\alpha$ is a prime power, then $\Omega$
tiles $\rr$ by translations.
\label{cor1}
\end{corollaire}

An easy application of Corollary \ref{cor1} yields the following.

\begin{corollaire}
Assume that $N=3$, then $\Omega$ is a spectral set if and
only if it tiles $\rr$ by translations.
\label{N=3}
\end{corollaire}

Note, however, that a spectrum need not always satisfy $\Lambda\subset
p^{-\alpha}\zz$, even if $N$ is a power of $p$. For instance, 
$\Omega= \{0,1,6,7\}+[0,1)$ has a spectrum $\Lambda=\{0,\frac{1}{12},
\frac{1}{2}, \frac{7}{12}\}+\zz$.

In light of the next theorem, the case when the degree of $A(x)$
is relatively small (equivalently, the set $\Omega$ is contained in a
relatively short interval) seems to be a natural  starting
point for further investigation.  The author expects that a more
general result should hold: if $\Omega\subset\rr^n$ ``almost fills up"
a convex set, and if it is a spectral set or tiles $\rr^n$
by translations, it must be a fundamental domain for a group.

\begin{theoreme}
Suppose that $A\subset[0,M-1]$, $M< 3N/2$, and let $\Omega=
A+[0,1)$.  Then the following are equivalent:

\smallskip
(i)  $A$ tiles $\zz$ by translations;

\smallskip
(ii) $\Omega$ is a spectral set;

\smallskip
(iii) $\Omega$ is a fundamental domain for the group $N\zz$,
\ie $A=\{0,1,\dots,N-1\}(\mod N)$;

\smallskip
(iv) $N\zz\oplus A=\zz$, \ie $N\zz\oplus \Omega=\rr$;

\smallskip
(v) $\Lambda=N^{-1}\zz$ is a spectrum for $\Omega$.
\label{thm-3N/2}
\end{theoreme}

The equivalence (iii) $\Leftrightarrow$ (iv) $\Leftrightarrow$
(v) (for an arbitrary domain $\Omega$) is Fuglede's theorem \cite{Fug},
and the implications (iv) $\Rightarrow$ (i), (v) $\Rightarrow$ (ii)
are trivial. Our new result is that the implications (i)
$\Rightarrow$ (iii) and (ii) $\Rightarrow$ (iii) hold for
$A$ as in the theorem.
Note that the set $\Omega_n=[0,n]\cup [2n,3n]$ tiles $\rr$
and is a spectral set, but is not a fundamental domain for any group;
this shows that the inequality $M<3N/2$ cannot be weakened.

We do not know whether a spectrum must always be rational, or
whether any spectral set must have a rational spectrum.
However, we have the following partial result.

\begin{theoreme}
Suppose that $\Omega=A+[0,1)$ has a spectrum $\Lambda$, and that
$\Omega\subset[0,M]$ for some $M<5N/2$. Then $\Lambda\subset\qq$.
\label{thm-2N}
\end{theoreme}

If the polynomial $A(x)$ is assumed to be irreducible, a stronger
result holds.

\begin{proposition}
(i) Suppose that $A$ tiles $\zz$ and that $A(x)$ is irreducible.
Then $N$ is prime, $A(x)=\Phi_{N^\alpha}(x)$ for
some $\alpha\in\nn$, $A=\{0,\alpha,2\alpha,\dots,(N-1)\alpha\}$,
and $\{k/N^{\alpha}:\ k=1,
2,\dots,N-1\}$ is a spectrum for $A(x)$.

\smallskip
(ii) Suppose that $A(x)$ is irreducible, has a spectrum, and
that deg$(A(x)) < \frac{5N}{2}$. Then the conclusions
of (i) hold, and, in particular, $A$ tiles $\zz$.
\label{prop-irred}
\end{proposition}

\medskip
If $m,n$ are integers, we will use $(m,n)$ to denote the
greatest common divisor of $m,n$. The cardinality of a finite
set $A$ is denoted by $\# A$. 

The author is grateful to Peter Borwein, David Boyd, and
Richard Froese for helpful conversations.

\section{Proof of Theorem \ref{CM-thm}}
\label{t1t2}
\init


We first prove (i).
Assume that (T1), (T2) hold, and define $S_A$ as in the
introduction.  
Consider the set $B$ of all numbers of the form:
\[
\sum_{s\in S_A} \frac{k_s}{s},
\]
where $k_s\in\{0,1,\dots,p-1\}$ if $s=p^\alpha\in S_A$ and $p$
is prime.  We claim that $B$ is a spectrum for $A$.

Let $N=p_1^{\aa_1}\dots p_m^{\alpha_m}$, where
$p_i$ are distinct primes. 
We have:
\beq
\Phi_s(1)=\left\{
\begin{array}{ll}
0\ &\hbox{ if }s=1,
\\
p\ &\hbox{ if }s=p^\alpha\hbox{ is a prime power,}
\\
1\ &\hbox{ otherwise }
\end{array}\right.
\label{eq.prime}
\eeq
(see \cite{CM}, Lemma 1.1). Thus (T1) implies that for each $i$
there are exactly $\alpha_i$ powers of $p_i$ in $S_A$. 
It follows that the cardinality of $B$ is
$p_1^{\aa_1}\dots p_m^{\alpha_m}=N.$

It remains to check that
\beq
A(e^{2\pi i(b-b')})=0\hbox{ for all }b,b'\in B.
\label{tt.claim}
\eeq
Fix $b,b'\in B$. Then $b-b'$ is a number of the form
$\sum_{s\in S_A} \frac{k_s}{s},$
where $k_{s}\in\{-p_i+1,-p_i+2,\dots,p_i-1\}$ if
$s$ is a power of $p_i$.  We rewrite this as
$b-b'=\sum_{i=1}^m b_i$,
where $b_i=0$ or
\[
b_i=\sum_{\beta:p_i^\beta\in S_A}\frac{k_{p_i^\beta}}{p_i^\beta},
\]
and $k_{p_i^\beta}\in\{-p_i+1,\dots,p_i-1\}$. If $b_i\neq 0$,
we may further write:
\[
b_i=\frac{k_i}{p_i^{\beta_i}},\ 
\]
where $\beta_i$ is the largest exponent such that $p_i^{\beta_i}\in S_A$,
$k_{p_i^{\beta_i}}\neq 0$, and $k_i\neq 0$ is an integer not divisible by $p_i$.
By Lemma \ref{tt.lemma1} below, 
\[
\Phi_{s}(e^{2\pi i(b-b')})=0\hbox{ for }s=\prod_{i:b_i\neq 0}p_i^{\beta_i}.
\]
This immediately implies (\ref{tt.claim}): since
$p_i^{\beta_i}\in S_A$ for each $i$, by (T2) $\Phi_s(x)$ divides $A(x)$,
hence $e^{2\pi i(b-b')}$ is also a root of $A(x)$.

\begin{lemme}
Let $b=\sum_{i=1}^n \frac{k_i}{s_i}$, where $(k_i,s_i)=1$ for all
$i$ and $(s_i,s_j)=1$ for all $i\neq j$.  Then 
$\Phi_s(e^{2\pi ib})=0$ for $s=s_1\dots s_n$.
\label{tt.lemma1}
\end{lemme}

\proof
The roots of $\Phi_s(x)$ are the primitive $s$-th roots of unity,
\ie the numbers $e^{2\pi ik/s}$ for $(k,s)=1$.  We write:
\[
b=\frac{k}{s},\ s=s_1\dots s_n,\ 
k=k_1\frac{s}{s_1}+\dots+k_n\frac{s}{s_n}.
\]
It suffices to check that $k,s$ are relatively prime. Suppose
therefore that there is a prime $p$ such that $p|k$ and $p|s$.
Then $p$ divides one of the $s_i$, say $s_1$. Since
$s_1$ divides $s/s_i$ for all $i\neq 1$, it follows that $p$ divides
$k_2\frac{s}{s_2}+\dots+k_n\frac{s}{s_n}$. Since we also assumed
that $p|k$, we must have $p|k_1\frac{s}{s_1}$, hence $p$ divides either
$k_1$ or $s/s_1$. This, however, contradicts the assumptions of
the lemma. Indeed, $p$ cannot divide $k_1$, since $(k_1,s_1)=1$. 
But we also have $(s_1,s_i)=1$ for $i\neq 1$, hence
$(s_1,s/s_1)=1$ and $p$ cannot divide $s/s_1$.
\qed

\bigskip
Next, we turn to part (ii) of the theorem.  Suppose that $\{\theta_j\}
\subset s^{-1}\zz$ is a spectrum for $A(x)$, where $s=p^\alpha$, 
$p$ is prime.  
Define $\theta_0=0$, $\theta_{ij}=\theta_i-\theta_j$,
then $e^{2\pi i\theta_{ij}}$ is a root of
$(x-1)A(x)$ for each $i,j=0,1,\dots,N-1$. 
We may write for all $i\neq j$:
\beq
\theta_{ij}=\frac{k_{ij}}{p^{\alpha_{ij}}},\ 
\alpha_{ij}\geq 1,\ (p,k_{ij})=1. 
\label{t.e1}
\eeq
Let $r_1<\dots<r_m$ be the distinct values of 
$\alpha_{ij}$ in (\ref{t.e1}). Since $\Phi_{p^{\alpha_{ij}}}(x)$
is the minimal polynomial of $e^{2\pi i\theta_{ij}}$, we must have
\[
\Phi_{p^{r_1}}(x)\dots\Phi_{p^{r_m}}(x)\ |\ A(x).
\]
From this and (\ref{eq.prime}) it follows that
\[
p^m=\Phi_{p^{r_1}}(1)\dots\Phi_{p^{r_m}}(1)\ |\ A(1)=N,
\]
and in particular $N\geq p^m$.  Thus (ii) will follow if we
prove that $N\leq p^m$. To do this, it suffices to verify that
the polynomial
\[
F_{0,r_1,\dots,r_m}(x)=
(x-1)\Phi_{p^{r_1}}(x)\dots\Phi_{p^{r_m}}(x)
\]
has at most $p^m$ roots $\eps_{i}$ such that 
\beq
F_{0,r_1,\dots,r_m}(\eps_i/\eps_j)=0\hbox{ for all }i,j.
\label{t.e2}
\eeq
The proof is by induction on $m$. If $m=1$, all roots of $F_{0,r_1}(x)
=(x-1)\Phi_{p^{r_1}}(x)$ are of the form $e^{2\pi ik/p^{r_1}}$,
$k\in\zz$.  If $k,k'$ belong to the same residue class (mod $p^{r_1-1}$),
then $k-k'$ and $p^{r_1}$ are not relatively prime, hence
either the corresponding roots $e^{2\pi ik/p^{r_1}}$ and 
$e^{2\pi ik'/p^{r_1}}$ are equal, or else $e^{2\pi i(k-k')/p^{r_1}}$
is not a root of $F_{0,r_1}(x)$. Since there are only $p$ distinct
residue classes (mod $p$), the claim follows for $m=1$.

Assume now that the claim is true with $m$ replaced by $m-1$. 
Let $\eps_j=e^{2\pi ik_j/r_m}$ be roots of $F_{0,r_1,\dots,r_m}(x)$
satisfying (\ref{t.e2}). We divide them into $p$ equivalence classes
$\Theta_0,\Theta_1,\dots,\Theta_{p-1}$:
\[
\eps_j=e^{2\pi ik_j/r^m}\in\Theta_l\Leftrightarrow
k_j=l(\mod p^{r_m-1}).
\]
By the same argument as above, if $\eps_i,\eps_j$ belong to the
same $\Theta_l$, $\eps_i/\eps_j$ cannot be a root of $\Phi_{p^{r_m}}(x)$,
and must therefore be a root of $F_{0,r_1,\dots,r_{m-1}}(x)=
(x-1)\Phi_{p^{r_1}}(x)\dots\Phi_{p^{r_{m-1}}}(x)$.
By the inductive hypothesis, $\#\Theta_l\leq p^{m-1}$.
Since the number of $\Theta_l$'s is $p$, the claim follows.
\qed


\section{Proof of Theorem \ref{thm-3N/2}}
\label{3N/2}
\init


Let $A=\{a_0=0,a_1,\dots,a_{N-1}\}\subset(\zz\cap[0,M-1])$,
$M<3N/2$, and let $\Omega=A+[0,1)$. We may assume that
$a_0<a_1<\dots<a_{N-1}$, $M=a_{N-1}+1$.  We must prove 
that each of the conditions (i), (ii) of Theorem 
\ref{thm-3N/2} implies that $\Omega$ is a fundamental domain
for $N\zz$, \ie $A=\{0,1,\dots,N-1\} (\mod N)$.

\bigskip
{\it Proof of (i) $\Rightarrow$ (iv).}
Assume that $A$ tiles $\zz$ by translations, \ie there
is a set $B=\{b_i\}_{i\in\zz}\subset \zz$ such that
$A\oplus B=\zz$ (and, consequently, $\Omega\oplus B=\rr$).
We may assume that $b_i<b_{i+1}$ for all $i$.
Define:
\[
c_i=b_i+M,\ \Omega_i=\Omega+b_i, \ I_i=[b_i,c_i).
\]
The interval $[b_i,c_{i+2})$ contains three disjoint translates of $\Omega$
($\Omega_i,\Omega_{i+1}, \Omega_{i+2}$), hence we must have
$c_{i+2}-b_i\geq 3N$. Using also that $M<3N/2$,
we see that:
\[
b_{i+2}-c_i=(c_{i+2}-M)-(b_i+M)\geq 3N-M>0.
\]
We also have the trivial inequality $c_i\geq b_{i+1}$ (if
it failed, $\Omega\oplus B$ would have ``gaps" $(c_i,b_{i+1})$.)
Hence each $I_i$ may overlap only with its immediate predecessor $I_{i-1}$
and successor $I_{i+1}$:
\beq
\dots <b_i\leq c_{i-1}<b_{i+1}\leq c_i<b_{i+2}\leq \dots.
\label{tt.e01}
\eeq

For $m>0$, define:
\[
S_1(m)=\sum_{i=1}^m (b_{i+1}-c_{i-1}),\ 
S_2(m)=\sum_{i=1}^m (c_i-b_{i+1}).
\]
Then:
\beq
\begin{array}{l}
\int_{b_1}^{c_m}\sum_{i=-\infty}^\infty\chi_{\Omega_i}
=S_1(m)+S_2(m)+O(1),
\\[3mm]
\int_{b_1}^{c_m}\sum_{i=-\infty}^\infty\chi_{I_i}
=S_1(m)+2S_2(m)+O(1),
\end{array}
\label{tt.e02}
\eeq
where we used that $\sum_{i=-\infty}^\infty\chi_{\Omega_i}\equiv 1$
and that, by (\ref{tt.e01}),
\[
\sum_{i=-\infty}^\infty\chi_{I_i}(x)=
\left\{
   \begin{array}{ll}
   2,\ &x\in (b_{i+1},c_i),
   \\
   1,\ &x\in (c_{i-1},b_{i+1}).
   \end{array}
\right.
\]
(Here and in the sequel, $O(1)$ denotes a quantity which, for
any fixed $N$ and $M$, is bounded uniformly in $m$ as $m\to\infty$.)
On the other hand, counting the $\Omega_i$'s contained in
$[b_1,c_m]$ and the $I_i$'s having non-empty intersection
with it, we obtain that:
\beq
\begin{array}{l}
\int_{b_1}^{c_m}\sum_{i=-\infty}^\infty\chi_{\Omega_i}
\geq mN+O(1),\\[3mm] 
\int_{b_1}^{c_m}\sum_{i=-\infty}^\infty\chi_{I_i}
\leq mM+O(1) \leq \alpha mN+O(1),
\end{array}
\label{tt.e03}
\eeq
for some $\alpha<3/2$ independent of $m$. From (\ref{tt.e02})
and (\ref{tt.e03}) we have:
\[
\begin{array}{l}
S_1(m)+2S_2(m)\leq\alpha mN+O(1),
\\[3mm]
S_1(m)+S_2(m)\geq mN+O(1),
\end{array}
\]
which yields that:
\[
S_2(m)\leq(\alpha-1)mN+O(1),\ 
S_1(m)\geq(1-\alpha)mN+O(1).
\]
In particular, $\lim_{m\to\infty}\frac{S_1(m)}{m}\geq (1-\alpha)N
>N/2$, hence we must have $b_{i+1}-c_{i-1}>N/2$ for some $i$.
But, by (\ref{tt.e01}), $(c_{i-1},b_{i+1})\subset \Omega_i\setminus
\bigcup_{ j\neq i}\Omega_j$, hence $\Omega$ contains at least one
``uninterrupted" interval of length $>N/2$.

Assume therefore that $[m,n]\subset\Omega$, $n-m>N/2$, and
that $[m,n]$ is maximal: $m-\half,n+\half\notin\Omega$.
Let $m_i=m+b_i$, $n_i=n+b_i$. By (\ref{tt.e01}),
\[
([m_i-1,m_i)\cup [n_i,n_i+1))\subset(I_{i-1}\cup I_{i+1}).
\]
Suppose that $[m'-1,m')\subset I_{i+1}$ for some $b_i\leq m'\leq m_i$.
Since $b_{i+1}>1$, $I_{i+1}$ contains at least one unit interval
to the right of $I_i$, and, in particular, $[m_i,n_i)\subset I_{i+1}$.
But $[m_i,n_i)\cap\Omega_{i+1}=0$. Thus $|I_{i+1}|\geq |\Omega_{i+1}|
+(n_i-m_i)>3N/2$, which contradicts the assumption that $M=|I_{i+1}|
<3N/2$.  We must therefore have $[m'-1,m')\subset I_{i-1}$. 
Similarly, $[n',n'+1)\subset I_{i+1}$ for all $n_i\leq n'\leq c_i-1$.
It follows that:
\[
c_{i-1}=m_i,\ b_{i+1}=n_i,\ 
[b_{i+1},c_i)\subset \Omega_i\cup\Omega_{i+1}.
\]
Hence:
\[
\begin{array}{l}
b_{i+1}-c_{i-1}=n_i-m_i=n-m,
\\[3mm]
c_i-b_{i+1}=
|\Omega_i\cap [n_i,c_i)|+|\Omega_{i+1}\cap[b_{i+1},c_i)|
=|\Omega\cap [n,M)|+|\Omega\cap[0,m)|
\\[3mm]
=|\Omega\setminus[m,n)|=N-(n-m),
\end{array}
\]
so that:
\[
c_{i}-c_{i-1}=(c_i-b_{i+1})+(b_{i+1}-c_{i-1})
=(N-(n-m))+(n-m)=N.
\]
This proves (iv).

\bigskip

{\it Proof of (ii) $\Rightarrow$ (v).} 
Suppose that $\Omega$ is a spectral set. By Proposition \ref{JP-prop},
$A(x)=\sum x^{a_i}$ has a spectrum $\{\theta_1,\dots,\theta_{N-1}\}
\subset (0,1)$. We will also denote $\theta_0=0$. Then:
\[
A(\eps_{ij})=0,\ i\neq j,
\]
where $\eps_{ij}=e^{2\pi i(\theta_i-\theta_j)}$.

We will first prove that
\beq
\forall i,j,k,l\  \exists m,i',l' \hbox{ such that }
\eps_{ij}=\eps_{i'm},\ \eps_{kl}=\eps_{ml'}.
\label{gp.e1}
\eeq
This will imply that $G=\{\eps_{ij}:\ i,j=0,1,\dots,N-1\}$
is a group. Indeed, since $G$ is finite and contains $1=\eps_{ii}$,
it suffices to prove that $G$ is closed under multiplication.
But if (\ref{gp.e1}) holds, then for any $\eps_{ij},\eps_{kl}\in G$
we have
\[
\eps_{ij}\eps_{kl}=\eps_{i'm}\eps_{ml'}=\eps_{i'l'}\in G.
\]

The verification of (\ref{gp.e1}) is as follows. Observe that:
\[
\eps_{ij}=\eps_{i'm}\Leftrightarrow \eps_{ii'}=\eps_{jm},\ 
\eps_{kl}=\eps_{ml'}\Leftrightarrow \eps_{ll'}=\eps_{km},
\]
hence it suffices to prove that:
\beq
\forall i,j,k,l\  \exists m,i',l' \hbox{ such that }
\eps_{ii'}=\eps_{jm},\ \eps_{ll'}=\eps_{km}.
\label{gp.e2}
\eeq
For each pair $(i,j)$ define:
\[
E_{ij}=\{m\in\{0,1,\dots,N-1\}:\ \eps_{jm}=\eps_{ii'}
\hbox{ for some }i'\},
\]
then (\ref{gp.e2}) is equivalent to:
\beq
E_{ij}\cap E_{kl}\neq\emptyset.
\label{gp.e3}
\eeq
Since $M<3N/2$, we have deg($A(x))<\frac{3N}{2}-1$. We have $\eps_{ii}=1$,
and $A(\eps_{ij})=0$ for all $i\neq j$, hence the number of the
possible values that $\eps_{ij}$ may take is bounded by 
\[
\hbox{deg}\,(A(x))+1<\frac{3N}{2}.
\]
On the other hand, since $\theta_i$ are distinct, so are
$\eps_{i0},\eps_{i1},\dots,\eps_{i,N-1}$ for any $i$. Hence for any
$i,j$, the number of values which occur both in
$\{\eps_{i0},\eps_{i1},\dots,\eps_{i,N-1}\}$ and in 
$\{\eps_{j0},\eps_{j1},\dots,\eps_{j,N-1}\}$ is $>N+N-\frac{3N}{2}
=N/2$. But this means that
\[
\# E_{ij}>N/2.
\]
Thus $E_{ij}, E_{kl}$ are two subsets of $\{0,1,\dots,N-1\}$
of cardinality $>N/2$, whereupon (\ref{gp.e3}) follows.

We have proved that $G$ is a group, hence $G$ must be the set of all
$n$-th roots of 1 for some $n\in\zz$. But $A(\eps_{ij})=0$ for all
$\eps_{ij}\in G$, hence the polynomial $P_n(x)=1+x+\dots+x^{n-1}$
divides $A(x)$, and, in particular, $n=P_n(1)|A(1)=N$. On the
other hand, since $\eps_i$ are all distinct, $n=\#G\geq N$.
Hence $n=N$, and the $\eps_{i}$ are $N$-th roots of 1, \ie 
$\Lambda=\frac{1}{N}\zz$.
\qed

\bigskip

Observe that the proof of the implication (ii) $\Rightarrow$ (v)
did not actually use that the degree of $A(x)$ was $<\frac{3N}{2}-1$,
but only that $A(x)$ has $<\frac{3N}{2}-1$ roots of modulus 1.
Hence we have also proved the following result, which
will be used in the proof of Theorem \ref{thm-2N}.

\begin{lemme}
Suppose that $A(x)$ has a spectrum, and that
$\#\{x: \ A(x)=0, |x|=1\}<\frac{3N}{2}-1$. Then
the conclusions (iii)--(v) of Theorem \ref{thm-3N/2} hold.
\label{roots-lemma}
\end{lemme}


\section{Proof of Theorem \ref{thm-2N}.}
\init


Suppose that $\Omega$ has a spectrum $\Lambda$, then
$A(x)$ has a spectrum $\{\theta_1,\dots,\theta_{n-1}\}$ as in
Proposition \ref{JP-prop}.  It suffices to prove that if
$\Omega\subset[0,M]$, and if $\theta_1$ is irrational, then
\beq
M\geq \frac{5N}{2}.
\label{m.e1}
\eeq
Suppose that $\theta_1\notin\qq$.  By Lemma \ref{roots-lemma},
$A(x)$ has at least $\frac{3N}{2}-1$ roots of modulus 1.  By
Lemma \ref{Kronecker} below, $A(x)$ has at least $N-1$ roots
off the unit circle.  Hence
\beq
M\geq\deg(A(x))+1\geq(N-1)+(\frac{3N}{2}-1)+1=\frac{5N}{2}-1.
\label{m.e2}
\eeq
We now improve this to (\ref{m.e1}) by a simple scaling argument.
Consider the set $k\Omega=kA=[0,k)$, where $k\in\nn$. Clearly,
if $\Omega$ has a spectrum $\Lambda\not\subset\qq$, then $k\Omega$
has a spectrum $k^{-1}\Lambda\not\subset\qq$.  But applying
(\ref{m.e2}) to $k\Omega$, we obtain that
\[
kM\geq \frac{5kN}{2}-1.
\]
Dividing both sides by $k$ and taking the limit $k\to\infty$,
we obtain (\ref{m.e1}).
\qed

\begin{lemme}
Suppose that $A(x)$ has a spectrum $\theta_1,\dots,\theta_{N-1}$,
and that $\theta_1\notin\qq$.  Then $A(x)$ has at least $N-1$ roots
of modulus $\neq 1$. 
\label{Kronecker}
\end{lemme}

\proof
Let $\eps_j=e^{2\pi i\theta_j}$, $\eps_0=1$, $\eps_{jk}:=\eps_j/\eps_k$;
then $\eps_j,j\neq 0$, and $\eps_{jk}, j\neq k$, are roots of $A(x)$.
Since $\eps_1$ is not a root of 1, by Kronecker's theorem \cite{Kron}
the minimal polynomial $A_1(x)$ of $\eps_1$ has at least one root
$\xi_1$ with $|\xi_1|\neq 1$.  

The Galois group $G$ of $A(x)$ acts transitively on the roots of
$A_1(x)$, hence there is a $\sigma\in G$ such that $\sigma(\eps_1) =\xi_1$.
Define $\xi_i=\sigma(\eps_i)$, $\xi_{ij}:=\xi_i/\xi_j=\sigma(\eps_{ij})$, 
then $\xi_j,j\neq 0$, and $\xi_{jk}, j\neq k$, are roots of $A(x)$.
Let
\[
I=\{i:\ 1\leq i\leq N-1,\ |\xi_i|=1\},\ k=\#I,
\]
and let $r_1,\dots,r_m$ be the distinct values of $|\xi_i|$ different
from $1$.  Since $|\xi_1|\neq 1$, we have $m\geq 1$.

For each value of $r_i$, fix a root $\xi_j$ with $|\xi_j|=r_i$, and
consider the roots
\beq
\xi_{jj'}:\ j'\in I.
\label{m.e3}
\eeq
This yields $k$ distinct roots of modulus $r_i$. The total number of 
roots (\ref{m.e3}), for all $r_i$, is $mk$.

We now consider two cases.

\begin{itemize}

\item If all $r_i$ are $>1$, consider the roots 
\[
1/\xi_i: \ i\neq 0, i\notin I,
\]
which have modulus $<1$ and are therefore distinct from the roots
found in (\ref{m.e3}).  The number of such roots is $N-1-k$.  Hence
the total number of roots of modulus $\neq 1$ is at least
$mk+(N-1-k)=N-1+(m-1)k\geq N-1$.  The case when all $r_i$ are $<1$ is
similar.

\item Suppose now that $r_2=\min (r_i)<1$, $r_3=\max(r_i)>1$, and fix
$\xi_2,\xi_3$ such that $|\xi_2|=r_2$, $|\xi_3|=r_3$.  Consider the roots
\beq
\{\xi_{2j}:\ |\xi_j|>1\}\cup \{\xi_{3j}:\ |\xi_j|<1\}.
\label{m.e4}
\eeq
The roots (\ref{m.e4}) have modulus either $<r_2$ or $>r_3$, hence
are distinct from those in (\ref{m.e3}). The number of roots in
(\ref{m.e4}) is $N-1-k$, so that again we have at least 
$mk+(N-1-k)\geq N-1$ roots of modulus $\neq 1$.
\qed

\end{itemize}


\section{Miscellaneous}
\label{sec.irred}
\init


{\bf Proof of Proposition \ref{prop-irred}.}
We begin with (i).
Suppose that $A(x)$ is irreducible, and that $A$ tiles $\zz$.
Then there are $B\subset\zz$, $m\in\zz$ such that $A(x)B(x)=
P_m(x)(\mod x^m-1)$, where $B(x)=\sum_{b\in B} x^b$ and $P_m(x)=
1+x+\dots+x^{m-1}$ (see \cite{CM}, Lemma 1.3). Since $A(x)$ is irreducible, 
either $A(x)$ divides $P_m(x)$ or $P_m(x)$ divides $B(x)$.
But in the latter case $m=P_m(1)=A(1)B(1)$ divides $B(1)$,
hence $N=A(1)=1$ and $A=\{0\}$.  Therefore we must have $A(x)|P_m(x)$, \ie
$A(x)=\Phi_k(x)$ for some $k|m$. By Theorem \ref{CM-thm},
\[
N=A(1)=\prod_{s\in S_A}\Phi_s(1)=\Phi_k(1).
\]
SInce $N\neq 1$, it follows from (\ref{eq.prime}) that $N$
is prime and $k=N^\alpha$ for some $\alpha\in\nn$. 
Then
\[
A(x)=1+x^\alpha+x^{2\alpha}+\dots+ x^{(n-1)\alpha},
\]
(cf. \cite{CM}, Lemma 1.1), and it is
trivial to verify that $\{0,\frac{1}{N^\alpha}, \frac{2}{N^\alpha},
\dots,\frac{N-1}{N^\alpha}\}$ is a spectrum for $A(x)$.

It remains to prove (ii). Suppose that $A(\eps)=0$ for at least one
root of unity $\eps$.  Since $A$ is irreducible, $A(x)=\Phi_s(x)$
for some $s$; using (\ref{eq.prime}) as above, we find that $N=A(1)$
is prime and hence the conclusions of (i) hold.  In particular,
$A=\{0,\alpha,2\alpha,\dots,(N-1)\alpha\}$ tiles $\zz$.

Assume therefore that $A(x)$ has a spectrum $\{\theta_1,\dots,\theta_{N-1}\}$,
and that $\theta_1$ is irrational. We shall first prove that $A(x)$ must
have at least $N$ roots off the circle $|x|=1$; combining this with
Lemma \ref{roots-lemma}, we find that 
\beq
\hbox{deg}\,(A(x))\geq \frac{5N}{2}-1.
\label{degree}
\eeq

Let $\eps_j=e^{2\pi i\theta_j}$.  Then $\eps_{j}$ and
$\eps_{jk}:=\eps_j/\eps_k$, $j\neq k$, are roots of $A(x)$,
and $\eps_1$ is not a root of unity. 
By Kronecker's theorem \cite{Kron}, $A(x)$
has at least one root $\xi$ with $|\xi|\neq 1$. Also,
$\bar\eps_1$ is a root of $A_1(x)$, hence by elementary Galois theory
so is $1/\xi$.  Thus $A(x)$ has at least one root of modulus
$>1$. 

Let $\xi_1$ be a root of $A(x)$ of maximal modulus; from the
previous paragraph we have $|\xi_1|>1$.  Since $A(x)$ is irreducible,
its Galois group $G$ is transitive, hence there is a $\sigma\in G$
such that $\sigma(\eps_1) =\xi_1$. Define $\xi_i=\sigma(\eps_i)$,
then $\xi_i$ and $\xi_{ij}:=\xi_i/\xi_j=\sigma(\eps_{ij})$ are roots of $A(x)$. 
Consider the sequence of roots:
\beq
\xi_1,\xi_2,\xi_3,\dots,\xi_{N-1}, \xi_{12},\xi_{13},\dots,\xi_{1,N-1}.
\label{irr.e1}
\eeq
Let $j\geq 2$. By the maximality of $|\xi_1|$, we have $|\xi_j|\leq|\xi_1|$
and $|\xi_{ij}|=|\frac{\xi_1}{\xi_j}|\leq|\xi_1|$, hence:
\[
1\leq |\xi_j|\leq |\xi_1|.
\]
Moreover, since $|\xi_1|\neq 1$, at most one of $\xi_j$,
$\xi_{1j}$ has modulus 1. Hence the sequence
\[
\xi_2,\xi_3,\dots,\xi_{N-1}, \xi_{12},\xi_{13},\dots,\xi_{1,N-1}.
\]
contains at least $N-2$ entries with modulus $>1$. 
Since the $\xi_i$ are distinct and $\neq 1$ for $i\neq 0$, the
value $\xi_1$ appears only once in (\ref{irr.e1}), and any other
value is taken at most twice (once as $\xi_j$ and once as $\xi_{1j'}$).
Hence (\ref{irr.e1}) contains at least $1+(N-2)/2=N/2$ distinct 
roots of modulus $>1$.
Similarly, using $1/\xi_1$ instead of $\xi_1$, we may find at least
$N/2$ distinct roots of modulus $<1$. This proves the claim that 
$A(x)$ has at least $N$ roots off the unit circle.

Finally, we improve (\ref{degree}) to 
\[
\deg(A(x))\geq \frac{5N}{2}.
\]
Observe that this improvement is automatic if $N$ is odd, since then
$N/2$ is not an integer and the above construction yields 
$\frac{N+1}{2}+\frac{N+1}{2}=N+1$ roots of modulus $\neq 1$.
Assume therefore that $N$ is even, and denote $|\xi_1|=R$. 
From the above construction, deg$(A(x))$ may equal $\frac{5N}{2}-1$
only if the following hold:

\medskip

(1) For any $2\leq i\leq N-1$, only one of $\xi_i$, $\xi_{1i}$
has modulus $\neq 1$. Hence $|\xi_i|=1$ or $R$.

\medskip

(2) For any $2\leq i\leq N-1$ there is a $2\leq j\leq N-1$ 
such that $\xi_1=\xi_{1j}$.  Since $\xi_i$ are distinct, this is
a 1-1 correspondence; moreover, (1) above implies that $i\neq j$.

\medskip
(3) There are no roots of $A(x)$ of modulus $R$ other than
those $\xi_i$ with $|\xi_i|=R$.

\medskip
By (1), (2), we must have $\{1,2,\dots,N-1\}=I_1\cup I_2$, where
$I_1=\{i: \ |\xi_i|=R\}$, $I_2=\{i: \ |\xi_i|=1\}$, 
$1\in I_1$, $\# I_1=N/2$, $\# I_2=\frac{N}{2}-1$.
Let $i\in I_1$, $j\in I_2$, then $|\xi_{ij}|=\frac{|\xi_i|}{
|\xi_j|}=R$, hence by (3) there is a $k=k(i,j)\in I_1$ such that
$\xi_{ij}=\xi_k$. Note that for any fixed $j$ the mapping
$i\to k(i,j)$ is one-to-one, and therefore may be inverted.
Thus, if $j\in I_2$ is fixed, for any $k\in I_1$ there
is a $i\in I_1$ such that $\xi_{ij}=\xi_k$, or, equivalently,
$\xi_k\xi_j=\xi_i$.  Iterating this procedure, we find that
for any $j\in I_2$ and $k\in I_1$:
\[
\xi_k,\ \xi_k\xi_j,\ \xi_k\xi_j^2,\dots, \xi_k\xi_j^n,\dots
\in\{\xi_i:\ i\in I_1\}.
\]
But $\#I_1=N/2$, hence $\xi_j^s=1$ for some $s\leq N/2$. 
Since $A(x)$ is irreducible, we must have $A(x)=\Phi_s(x)$ for
some $s\leq N/2$, which clearly contradicts our assumptions.
\qed

\bigskip


{\bf Proof of Corollary \ref{N=3}.} Let $A=\{0,a_1,a_2\}$. We may assume
that $(a_1,a_2)=1$: indeed, let $k\in\zz$, then $A$ tiles $\zz$ if and only 
if $kA$ tiles $\zz$ (\cite{CM}, Lemma 1.4(1)), and $\{\theta_1,\theta_2\}$
is a spectrum for $A(x)$ if and only if $\{\theta_1/k,\theta_2/k\}$ is
a spectrum for $(kA)(x)=A(x^k)$. 

If $A$ tiles $\zz$, $A(x)$ has a spectrum by Corollary \ref{cor1}(i).
Conversely, suppose that $A(x)$ has a spectrum $\{\theta_1,\theta_2\}$, then:
\[
1+e^{2\pi ia_1\theta_j}+e^{2\pi ia_2\theta_j}=0,\ j=1,2.
\]
Hence $e^{2\pi ia_1\theta_j}$, $e^{2\pi ia_2\theta_j}$ are cubic
roots of $1$, and in particular $a_1\theta_j, a_2\theta_j\in\frac{1}{3}\zz$.
Since $(a_1,a_2)=1$, there are integers $k_1,k_2$ such that $k_1a_1
+k_2a_2=1$, so that
\[
\theta_j=k_1a_1\theta_j +k_2a_2\theta_j\in\frac{1}{3}\zz.
\]
By Corollary \ref{cor1}(ii), $A$ tiles $\zz$.
\qed


\end{document}